\documentclass{amsart}

\usepackage{amsthm}

\hypersetup{
  colorlinks=true,
  urlcolor=blue,
  linkcolor=blue,
  citecolor=blue}

\newtheorem{theorem}{Theorem}[section]

\theoremstyle{definition}

\theoremstyle{remark}

\numberwithin{equation}{section}

\begin{document}

\title{Varieties of mathematical understanding}

\author{Jeremy Avigad}
\address{Department of Philosophy, Baker Hall 161, Carnegie Mellon University, Pittsburgh, PA 15213, USA}
\curraddr{}
\email{avigad@cmu.edu}
\thanks{I am grateful to Robert Y.~Lewis, Patrick Massot, and the anonymous referee for comments, corrections, and helpful suggestions. This work has been partially supported by AFOSR grant FA9550-18-1-0120 and the Sloan Foundation.}


\date{February 21, 2021}

\dedicatory{}

\begin{abstract}
This essay considers ways that recent uses of computers in mathematics challenge contemporary views on the nature of mathematical understanding. It also puts these challenges in a historical perspective and offers speculation as to a possible resolution.
\end{abstract}

\maketitle

\section{Introduction}

The goal of this essay is to explore some of the ways that recent uses of computers in mathematics challenge conventional views as to the nature of mathematical understanding. Section~\ref{section:theorems} surveys four developments in discrete geometry, all of which rely on the use of computation, although in markedly different ways. Section~\ref{section:understanding} aims to clarify the view that an important goal of mathematics is to obtain a deep \emph{conceptual} understanding. Section~\ref{section:assessment} highlights some of the differences between such a conceptual understanding and the kinds of computational understanding we get from the results in discrete geometry. Finally, Section~\ref{section:discussion} reflects on the history of mathematics to support speculation as to whether some sort of reconciliation is possible, and, if so, what form it might take.

\section{Some results in discrete geometry}
\label{section:theorems}

Discrete geometry studies combinatorial properties of discrete geometric objects. A number of important problems in the field have to do with optimal or otherwise interesting arrangements of such objects in space. Here, ``space'' might mean the Euclidean plane or three-dimensional Euclidean space, but it can also refer to higher-dimensional and non-Euclidean spaces. In this section, we will consider four illustrative developments, each centered on a representative theorem.

\subsection{The Kepler conjecture}

The first result we will consider was known for a long time as the \emph{Kepler conjecture}, but now perhaps should go by the name \emph{Hales' theorem}. Imagine filling three-dimensional space with infinitely many nonoverlapping spheres, all the same size. One can arrange them in hexagonal layers, with each successive layer nestled into the one beneath it. At each layer, there are two choices of how to fit the one above it into the crevices, but they are equally efficient; the \emph{face-centered cubic packing} is one such arrangement. It is not hard to show that the fraction of space filled in this way is $\pi / 3 \sqrt 2$, about $0.7405$.

In 1611, Johannes Kepler asserted that this density is optimal, but, for centuries, the claim went unproved. In 1998, Thomas Hales announced a proof, obtained with a contribution from a student of his, Samuel Ferguson \cite{hales:05}. Hales' strategy was to reduce the problem to a complex set of sufficient conditions that could be verified by an extensive body of computation. The proof shows that any counterexample to the Kepler conjecture would imply the existence of a finite arrangement of spheres satisfying certain properties, giving rise, in turn, to a certain combinatorial structure. Computer code then produced an exhaustive enumeration of the possible combinatorial structures; to be realized geometrically, any such structure would have to satisfy certain inequalities. Using branch-and-bound methods, these inequalities were relaxed to linear ones, at which point linear programming methods were used to demonstrate their infeasibility. In other words, the computations showed that there is no finite arrangement of spheres of the kind guaranteed by a putative counterexample. The proof thus consisted of a traditional mathematical argument (about 250 pages at the time) combined with a substantial body of computer code used to carry out the calculations.

\begin{theorem}
\label{theorem:hales}
The optimal density achieved by a packing of nonoverlapping equally sized spheres in Euclidean space is the one attained by the face-centered cubic packing.
\end{theorem}
\noindent Hales' \emph{Flyspeck} project, completed in 2014, used computational proof assistants to verify the correctness of the result in a formal axiomatic framework \cite{hales:et:al:17}.

\subsection{Packings of regular tetrahedra}

The second theorem we will consider has to do with the optimal packing density of regular tetrahedra instead of spheres. In 2006, John Conway and Salvatore Torquato obtained a packing with a density of about $0.7175$, and their paper sparked a flurry of interest in the problem \cite{conway:torquato:06}. The best result to date is due to Elizabeth Chen, Michael Engel, and Sharon Glotzer \cite{chen:engel:glotzer:10}:
\begin{theorem}
\label{theorem:tetrahedra}
There is a packing of equally sized regular tetrahedra in Euclidean space with density $4000/4671$.
\end{theorem}
\noindent That fraction is equal to about $0.8563$.

Theorem~\ref{theorem:tetrahedra} emerged from an impressive constellation of approaches to the problem. Conway and Torquato came up with their packing by starting with a configuration of irregular tetrahedra that fills space exactly, approximating it with regular tetrahedra, and then tweaking the configuration to increase the density. Other researchers adopted the strategy of parameterizing families of packings that seemed promising. For example, one can focus on packings based on a two-tetrahedron \emph{dimer}the configuration formed by gluing two tetrahedra together so that their bases match. With such restrictions, researchers could use various computer search techniques to optimize the choice of parameters.

Yet a third approach used Monte Carlo simulations of roughly the following physical act: put a lot of tetrahedra into a big box, jiggle the box, and then compress it as the configuration settles. This often results in arrangements that exhibit a surprisingly orderly quasi-crystalline structure. There is a thermodynamic explanation: because pockets of order are locally more efficient, their components have more freedom to move around, resulting in a system with higher entropy overall. (Intuitively, there are lots of ways to perturb the components, so such arrangements dominate the space of configurations, and a random process is more likely to end up there.) In 2009, Glotzer and colleagues used such a simulation (and subsequent tuning) to find a packing with a density of $0.8503$, a record at the time \cite{haji:akbari:et:al:09}. The current record of $0.8563$ was found using a parameterized family instead, but similar experimental techniques were used to increase faith in that the result is optimal at least with respect to nearby alternatives.

\subsection{Sphere packing in higher dimensions}

Packing problems are not restricted to three dimensions, and the next results have to do with the densest packing of spheres in $n$-dimensional Euclidean space. Surprisingly little is known about dimensions greater than three, but not for lack of effort. There are remarkably efficient packings in eight and twenty-four dimensions. The eight dimensional packing centers spheres on the points of a lattice known as $E_8$, while the twenty-four dimensional packing centers spheres on the points of a lattice known as the Leech lattice, discovered by Conway in 1967. In 2016, Maryna Viazovska announced a proof that the density achieved by the $E_8$ packing is optimal in eight dimensions \cite{viazovska:17}, and within a week of the announcement, Henry Cohn, Abhinav Kumar, Stephen D.~Miller, and Danylo Radchenko had joined her to extend the method to show the optimality of the packing based on the Leech lattice in twenty-four dimensions \cite{cohn:et:al:17}.

\begin{theorem}
\label{theorem:leech}
The optimal density of an 8-dimensional sphere packing is attained by $E_8$, and the optimal density of a 24-dimensional sphere packing is attained by the Leech lattice.
\end{theorem}

There is a nice exposition of these results and the methods behind them in an article by Cohn in the \emph{Notices of the American Mathematical Society} \cite{cohn:17}. In a paper in the \emph{Annals of Mathematics} in 2003, Cohn and Noam Elkies provided a simple Fourier-analytic method of providing upper bounds on the density of sphere packings \cite{cohn:17}. In short, they showed that the existence of a function with certain properties gives rise to such a bound, thereby reducing the task to that of finding functions with the right properties. Using numerical methods, they were able to produce functions yielding the best known bounds for dimensions 4 through 36, and noted that the upper bounds were especially tight in dimensions 8 and 24. This gave rise to the hope that one could construct suitable ``magic functions'' in those dimensions to show that the lower bounds given by $E_8$ and the Leech lattice are sharp.

As the \emph{Notices} article explains, numeric calculations by Cohn and Miller provided hints as to what such functions should look like. Viazovska's construction for $E_8$ combined these insights with mathematical intuition and expertise: drawing on modular forms, properties of the Laplace transform, experimentation, and guesswork, she was able to produce a function with the requisite properties. The extension to 24 dimensions was based on the same techniques.

\subsection{The Keller conjecture}

The last result that we will consider starts with two dimensions and moves to higher dimensions from there. It is clearly not possible to tile the plane with unit squares in such a way that no two squares share a common edge. Start with a square and then stack one on top of it, shifted to the right, say, so that they will not have an edge in common. Put another square to the right of the first one, shifted down to avoid sharing an edge. There is then a gap left over that clearly cannot be filled. It is not hard to convince oneself that, similarly, there is no way to fill space with unit cubes without two cubes sharing a face. In 1930, Ott-Heinrich Keller conjectured that the corresponding claim holds in all dimensions, that is, for every $n$, there is no tiling of $n$-dimensional space with unit cubes such that no two squares share a common $(n-1)$-dimensional face \cite{keller:30}.

In 1940, Oskar Perron showed that the conjecture is true for $n \le 6$ \cite{perron:40a,perron:40b}. In 1992, Jeffrey Lagarias and Peter Shor showed, surprisingly, that it is false for $n \ge 10$ \cite{lagarias:shor:92}. (The failure of the conjecture in any dimension straightforwardly implies the failure in higher dimensions.) In 2002, John Mackey showed the conjecture is false for $n \ge 8$ \cite{mackey:02}, leaving only the case $n = 7$ open. This was recently settled by Joshua Brakensiek, Marijn Heule, John Mackey, and David Narv\'aez \cite{brakensiek:et:al:20}.

\begin{theorem}
\label{theorem:keller}
Keller's conjecture is true up to and including dimension 7, but fails for dimension 8 and above.
\end{theorem}

Once again, the final result builds on a heterogeneous mix of methods. An important milestone in the attack on the conjecture was to reduce the original geometric question to a discrete, graph-theoretic one. In 1986 S\'andor Szab\'o reduced the conjecture to periodic tilings \cite{szabo:86}, and in 1990 Kereszt\'ely Corr\'adi and Szab\'o showed that the conjecture is equivalent to the statement that for every $n$, there is no clique of size $2^n$ in a certain finite graph, $\Gamma_n$. In fact, that was how the initial counterexamples were obtained: Lagarias and Shor exhibited a clique of size $2^{10}$ in $\Gamma_{10}$, and Mackey exhibited a clique of size $2^8$ in $\Gamma_{8}$, each of which gave rise to a counterexample of the corresponding dimension. In 2017, Andrzej Kisielewicz generalized the construction to a family of graphs $G_{n, s}$ to make the correspondence with the dimension sharp: the Keller conjecture is true at dimension $n$ if and only if there is no clique in $G_{n, s}$ for any $s$. (In this notation, the Corr\'adi--Szab\'o graph $\Gamma_n$ is $G_{n, 2}$.) He then showed that, to settle the only remaining case, it is enough to rule out the existence of a clique in $G_{7, 6}$, which has $12^7$ nodes \cite{kisielewicz:17}.

That is exactly what Brakensiek, Heule, Mackey, and Narv\'aez did. To obtain the result, they used propositional logic to encode the existence of a clique in $G_{7,6}$ and then used software known as a \emph{SAT solver} to show that the resulting formula has no solutions. Propositional satisfiability solvers have, in recent years, opened up new possibilities for finding complex combinatorial configurations meeting various sets of constraints, or showing that a given set of constraints cannot be met \cite{biere:et:al:09}. They are used to solve a wide range of problems in hardware verification, software verification, and AI. A naive search would have overwhelmed the solver; to obtain the solution, Brakensiek, Heule, Mackey, and Narv\'aez had to develop new heuristics to break symmetry and reduce the search space to a manageable size.

\section{Conceptual understanding and depth}
\label{section:understanding}

What do results like these contribute to mathematical understanding? Answering that question requires a conception of the kind of understanding that mathematics is after. The subject comprises a range of pursuits, some more pure, and some more applied; it reaches into all the physical and social sciences, including engineering, statistics, and computer science. What I am trying to get at here is something like the essence of pure mathematics, the things we value about mathematics itself. What makes something worthy of publication in a premier journal like the \emph{Annals of Mathematics}? What are the criteria at play in selecting a Fields medalist? What are the topics most worthy of representation at the International Congress of Mathematicians? In short, what is quintessentially important to the subject?

It is common to say that the goal of mathematics is to obtain a \emph{conceptual} understanding of mathematical phenomena, and a \emph{deep} understanding at that. The Langlands program, which seeks to develop far-reaching connections between number theory and geometry, is often held as a paradigm of conceptual depth. What makes it so? To develop some intuitions, I will draw on informal writings by Kevin Buzzard, a number theorist at Imperial College in London. In 2017, Buzzard launched his \emph{Xena} blog,\footnote{\url{https://xenaproject.wordpress.com/}} in part to document his newfound interest in the use of computational proof assistants to develop libraries of formally checked mathematics. Over the years, he has used it as a venue to share what is on his mind, sometimes explaining to computer scientists and the general public what is important to him as a mathematician, and sometimes explaining to mathematicians what he has learned about proof assistants. The first of these two activities is the one that is relevant here. One does not have to look far to find other mathematicians expressing similar values, but Buzzard has a gift for exposition, and it is especially convenient to find so many thoughtful characterizations of mathematical activity collected in one place.\footnote{Philosophers of mathematics have weighed in as well. See, for example, the special issue of \emph{Philosophia Mathematica} devoted to the notion of mathematical depth (volume 23, issue 2, June 2015), or a recent survey by Hamami and Morris \cite{hamami:morris:20}, which aims to provide mathematics educators with an overview of recent philosophical work on mathematical representation, explanation, understanding, and design.}

One observation is that deep mathematics is usually pretty complicated. Generally speaking, the deeper the result, the harder it is for the general public or even mathematicians not directly involved with the research to appreciate it.
\begin{quote}
  So what are the mathematicians I know interested in? Well, let’s take the research staff in my department at Imperial College. They are working on results about objects which in some cases take hundreds of axioms to define, or are even more complicated: sometimes even the definitions of the objects we study can only be formalised once one has proved hard theorems. For example the definition of the canonical model of a Shimura variety over a number field can only be made once one has proved most of the theorems in Deligne's paper on canonical models, which in turn rely on the theory of CM abelian varieties, which in turn rely on the theorems of global class field theory. That's the kind of definitions which mathematicians in my department get excited about [\ldots]. I once went to an entire 24 lecture course by John Coates which assumed local class field theory and deduced the theorems of global class field theory. I have read enough of the book by Shimura and Taniyama on CM abelian varieties to know what's going on there. I have been to a study group on Deligne's paper on canonical models. So after perhaps 100 hours of study absorbing the prerequisites, I was ready for the definition of a Shimura variety over a number field. And then there is still the small matter of the definition of \'etale cohomology. (Xena, July 6, 2018)
\end{quote}
Perhaps more important than the complexity of statements and definitions is the complexity of the proofs. A really deep proof often requires more background knowledge than any one person can master.
\begin{quote}
  To completely understand a proof of FLT (let’s say, for now, the proof explained in the 1995 Darmon-Diamond-Taylor paper) you will need to be a master of the techniques used by Langlands in his proof of cyclic base change (and I know people who are), and a master of Mazur's work on the Eisenstein ideal (and I know people who are). But what about the far less sexy technical stuff? To move from the complex analytic theory of modular forms to the algebraic theory you will need to know the delicate interplay between the analytic and algebraic theory of moduli spaces of elliptic curves (and I know people who know this---but I went through some of this stuff once and it's far more delicate than I had imagined, and there are parts of it where the only useful reference seems to be Brian Conrad's brain). This last example is perhaps a good example of a tedious technical issue which it’s very easy to forget about, because the results are intuitive and the proofs can be technical. There are many other subtleties which one would have to fully understand because they're on the syllabus. Is there really one human being who would feel confident answering questions on all of this material? I am really not sure at all. (Xena, September 27, 2019)
\end{quote}

It would be a mistake, however, to equate depth with complexity, and other postings on Xena make it clear that the complexity is only a means to an end. Complex definitions and proofs are worth the effort when they provide answers to questions that are judged by the community to be interesting and important. This last standard is one that is hard to pin down, but, as Buzzard emphasizes, that does not mean that it is arbitrary or subjective.

\begin{quote}
\noindent What mathematics is fashionable? Just take a look at the work of the recent Fields Medalists. That's a pretty good way of telling.

But fortunately, unlike many other fashions, ``fashionable mathematics'' is not controlled by the whim of big companies or some cabal. Fashionable mathematics is mathematics which justifies itself by its ability to answer questions which have been previously deemed interesting or important. Peter Scholze’s definition of a perfectoid space opened a door. In the last ten years, perfectoid spaces have been used to prove new cases of the monodromy--weight conjecture, to prove the direct summand conjecture, to give a new proof of purity for flat cohomology, a strengthened version of the almost purity theorem and so on (and I didn't even mention applications to the Langlands philosophy). These are results whose statements do not mention perfectoid spaces, and some have been open for decades. This is what makes this kind of mathematics fashionable---it is giving my community new insights. (Xena, February 9, 2020)
\end{quote}

Buzzard's posts highlight the importance of having the right definitions, as does Scholze himself in a lovely quotation featured on Michael Harris' \emph{Mathematics without Apologies} blog:\footnote{\url{https://mathematicswithoutapologies.wordpress.com}}
\begin{quote}
What I care most about are definitions. For one thing, humans describe mathematics through language, and, as always, we need sharp words in order to articulate our ideas clearly. (For example, for a long time, I had some idea of the concept of diamonds. But only when I came up with a good name could I really start to think about it, let alone communicate it to others. Finding the name took several months (or even a year?). Then it took another two or three years to finally write down the correct definition (among many close variants). The essential difficulty in writing ``\'Etale cohomology of diamonds'' was (by far) not giving the proofs, but finding the definitions.) But even beyond mere language, we perceive mathematical nature through the lenses given by definitions, and it is critical that the definitions put the essential points into focus. (Mathematicians without Apologies, June 2, 2018)
\end{quote}
The last sentence shows that it isn't even the definitions per se, but the intuitions they convey. Sometimes this is cast in terms of providing new representations:
\begin{quote}
  Mathematicians are so good at instantly switching between the various ``obviously equivalent'' ways that a mathematician looks at a complicated algebraic object. (``It’s an equivalence relation! Now it’s a partition! Now it’s an equivalence relation again! Let your mental model jump freely to the point of view which makes what I'm saying in this particular paragraph obvious!'', or ``Matrices are obviously associative under multiplication because functions are associative under composition.'') (Xena, July 23, 2020)
\end{quote}
More generally, the value of a good definition is that it serves as the carrier for mathematical insight and ideas.
\begin{quotation}
  \noindent Ideas are the purely artistic part of mathematics. That ``lightbulb'' moment, the insight which enables you to solve a problem---this is the elusive mathematical idea.

  Ideas are the part of mathematics that I understand the least, in a formal sense. Here are two questions:
  \begin{itemize}
  \item What is a group?
  \item How do you think about groups?
  \end{itemize}
  The first one is a precise ``scientific'' question. A group is a set equipped with some extra structure, and which satisfies some axioms. The formal answer is on Wikipedia's page on groups. A group is a definition. But the second question is a different kind of question. Different people think about groups in different ways. Say $G$ is a group generated by an element $x$ satisfying $x^5=x^8=1$. What can you say about $G$? If you are a mathematics undergraduate who has just seen the formal definition of a group, you can probably say nothing. If you have a more mature understanding of group theory, you instantly know that this group is trivial, because you have a far more sophisticated model of what is going on. Ideas are complicated, and human-dependent. A computer's idea of what a group is, is literally a copy of the definition in Wikipedia, and this is one of the reasons that computers are currently bad at proving new theorems by themselves. You can develop a computer's intuition by teaching it theorems about groups, or teaching it examples of groups, or trying to write AI's which figure out group theory theorems or examples of groups automatically. But intuition is a very subtle thing, and I do not understand it at all well\ldots\,. (Xena, June 20, 2020)
\end{quotation}

If ideas and representations are the message, the language of perfectoid spaces and Shimura varieties is the medium. Conceptual understanding is couched in a structural world-view that is a hallmark of modern mathematics. The approach can be characterized in terms of all of the following: a focus on axiomatically described structures; the use of certain tools, such as limit and quotient constructions, for building new structures from old ones; and an emphasis on characterizing structures and their properties in relation to other structures, for example, by studying the morphisms between them and by viewing the structures themselves as elements of even more abstract structures.

One can discern the first glimmerings of the conceptual view of mathematics in Gauss' monumental \emph{Disquisitiones Arithmeticae} of 1801. That work appeared well before the advent of modern algebra---for example, Cayley's abstract definition of a group did not appear until 1854---but Gauss' second proof of the law of quadratic reciprocity and his methods of classifying binary quadratic forms hint at such abstractions. Early in the \emph{Disquisitiones} Gauss provided a fairly modern treatment of Wilson's theorem, which states that the product of the nonzero residues modulo a prime number $p$ is congruent to $-1$ modulo $p$. He attributed Wilson's inability to prove the theorem to an excessive reliance on notation, and wrote: ``in our opinion truths of this kind should be drawn from the notions involved rather than from notations.''\footnote{Clark's translation \cite[Section 76]{gauss:disquisitiones:arithmeticae} uses ``ideas'' instead of ``notions,'' sacrificing the play on words that is apparent in the original Latin: ``At nostro quidem iudicio huiusmodi veritates ex notionibus potius quam ex notationibus hauriri debebant.''} The slogans ``notions over notation'' and ``concepts over calculation'' were rallying cries for the nineteenth-century transformation of mathematics. Much as been written about the tensions between the computational rigor of Kronecker, Weierstrass, and the Berlin school, on the one hand, and more algebraic, geometric, set-theoretic, and conceptual approaches developed by Riemann, Dedekind, Cantor, and Hilbert on the other.\footnote{Dedekind quoted Gauss' slogan approvingly in 1895, in one of his later works on the theory of ideals. For an analysis of how ``concepts over calculation'' played out in Dedekind's work in algebraic number theory, see \cite{avigad:06}.} The conceptual movement grew to maturity in the nineteenth century, with the work of Emmy Noether, the publication of van der Waerden's \emph{Modern Algebra} in 1930, the founding of Bourbaki, and the work of Grothendieck. The transformation of mathematics has been so thorough that few mathematicians will even notice the extent to which the language of structures pervades Buzzard's posts. It is simply the way we talk about mathematics today.

Let me summarize some of the morals we can extract from Xena. What are judged to be the important questions facing mathematics change over time, but the changes are not random or arbitrary. Mathematics may have been designed to help us make sense of our physical experiences, but doing so gives rise to new mathematical experiences, which in turn give rise to new questions and new research programs. The fact that questions are situated in a tradition that extends through the centuries is important to judgments of depth. So is complexity: questions often turn out to be surprisingly hard, and when complex definitions and proofs are deemed to be necessary to address them, their complexity is a further sign of depth. The central challenge is to find the right definitions, which provide new representations, new formulations, and new contextualizations of previous knowledge. The ability of a research program to draw on disparate intuitions is significant, as is the ability to unify previously existing fields of study and open up fruitful new avenues of research. Finally, social aspects of mathematics are important: the fact that some results are too complex for any one person to master is evidence that we are pushing against the boundaries of what can be known, and part of what we value about mathematics is that it supports a collaborative understanding.

It may be helpful to think of conceptual mathematics as a technology that is designed to extend our cognitive reach. Nineteenth-century algebraic methods emerged to manage complexity in situations where calculational approaches had become unwieldy. Structural methods provide a language for defining and reasoning about abstract objects, which in turn, provide an effective means of controlling information by modularizing arguments, making information salient when it is needed, and suppressing it when it is a distraction. Depth reflects the extent to which our conceptual engineering succeeds in that regard, enabling us to reason better and more efficiently. Depth can also be viewed as a counterweight to cleverness: when it comes to solving really hard problems, cleverness can only get us so far, and conceptual infrastructure and scaffolding are often needed to get the job done. To be sure, progress often requires clever tactical skirmishes, but far-reaching conceptual strategy provides the most substantial means of extending our knowledge.

We can now assess the Langlands program in these terms. The program has served to unify several previously existing fields of mathematical research, including local and global class field theory, the classical theory of modular forms and Eisenstein series, the theory of L-functions, harmonic analysis on reductive groups, and Shimura varieties. It has contributed to major results in number theory, including Weil's conjecture on Tamagawa numbers, the modularity theorem, Fermat's last theorem, and the Sato-Tate conjecture. With the associated Langlands conjectures, it provides an ambitious but concrete research program, and it has inspired major imitations, such as the geometric Langlands program and ``mod $p$'' variants.\footnote{I am grateful to the anonymous referee for the contents of this paragraph.}

The program exhibits all the characteristics that are underscored by Xena. It has led to the resolution of existing problems, it supports the transfer of ideas between different fields of inquiry, and it has opened up new questions and new lines of research. It is historically linked to the geometric and number-theoretic traditions that stretch back to Euclid. It requires building massively elaborate networks of structures, with definitions that are often too complex to be written down in complete detail, and proofs that are even more so. The program draws on more mathematics than any one mathematician can master, but with modularization of concerns and the coordinated efforts, it enables us to make progress as a community. These are all indications of depth.

\section{Assessing discrete geometry}
\label{section:assessment}

How do the theorems of Section~\ref{section:theorems} fare with respect to the values described in Section~\ref{section:understanding}? Not well, I am afraid, but it is illuminating to see why not. It is clear from the start that the results fail to exhibit many of the features that are praised in Buzzard's posts. The definitions involved are not very complex and the proofs are fairly self-contained; the average mathematician can generally come up to speed with a week or two's worth of background reading and a handful of prior journal articles.

The results are further handicapped by virtue of their subject matter. Discrete geometry is closely allied with combinatorics, a subject which, in the eyes of many mathematicians, is a repository of cleverness rather than depth. William T.~Gowers has observed that, in the traditional division of mathematics into the two cultures of ``theory builders'' and ``problem solvers,'' the former are generally held in higher esteem \cite{gowers:00}. He has eloquently defended combinatorics against the criticism that it is an ad hoc compilation of clever but disjointed results, the result of beavering away at insignificant problems rather than contributing to the grand mathematical edifice. The status of combinatorics was raised by a number of impressive results after the turn of the new millennium, including the celebrated Green--Tao result on primes in an arithmetic progression \cite{green:tao:08}, as well as generalizations and strengthenings of Szemer\'edi's theorem. Such work drew on ideas and methods from ergodic theory and Fourier analysis, bringing combinatorics closer to the status of structural mathematics, and making it, in Buzzard's terminology, fashionable.\footnote{A use of the Hahn-Banach theorem by Gowers in the finite setting is a nice example of a reconciliation of infinitary and finitary methods \cite{gowers:10}.} But while this gave combinatorics a welcome boost, discrete geometry did not rise with it; the latter still looks like combinatorics used to look before its structural makeover.

But the results of Section~\ref{section:theorems} are even more tainted by their heavy use of computation, including optimization techniques, combinatorial enumeration, validated numerical computation, linear programming methods, Monte Carlo simulation, search techniques, propositional satisfiability algorithms, and computer algebra. These are mainstays of computer science, engineering, and applied mathematics, but they are not seen as the purview of pure mathematics. So it hard to view the results of Section~\ref{section:theorems} as part of the pure tradition.

When it comes to assessing the theorems of Section~\ref{section:theorems} with respect to the core values of understanding and depth, Theorem~\ref{theorem:tetrahedra}, the lower bound on the optimal density of tetrahedral packings, is the easiest to dismiss. Although the results required insight and hard work, there is not much theory there---no complex definitions and proofs, and no overarching conceptual framework. The result isn't known to be optimal, and there seems to be nothing special about the fraction $4000/4671$. And some of the related results, for example, the fact that certain types of random processes settle into certain types of ordered configurations, are only experimentally confirmed, though the claims are also supported by general physical principles. I once heard Sharon Glotzer present an engaging survey of the work at a conference to celebrate Hales' sixtieth birthday. At one point in her talk, after stating a result, she said ``we can show that rigorously,'' and then paused. Looking out at an audience of discrete geometers, number theorists, representation theorists, computer scientists, and logicians, she added ``well, rigorously for us, not for you.'' This tongue-in-cheek remark was a gracious way of acknowledging methodological differences between those fields and her own, chemical engineering.

Results like Theorem~\ref{theorem:tetrahedra} are sometimes classified as \emph{experimental mathematics}, a term which describes a grab-bag of ways of using computation to explore mathematical phenomena and discover new results. Sometimes results that are discovered using such techniques can be justified rigorously, after the fact; sometimes the computations are meant only to indicate patterns or lend plausibility to claims. The field, surveyed in a number of accessible textbooks \cite{borwein:bailey:08,borwein:devlin:09}, is lively and popular, and there is even a journal, \emph{Experimental Mathematics}, devoted to it. But to the mathematical elite, the phrase ``experimental mathematics'' is itself a contradiction in terms. Many mathematicians look down on such work as being recreational at best, and, at worst, detrimental to the subject.

Theorem~\ref{theorem:keller}, the resolution of the Keller conjecture, is the next easiest to dismiss. It has a few more things going in its favor: the conjecture had been around since 1930, algebraic insights were used to reduce the problem to a combinatorial one, and the nature of the results---the fact that something counterintuitive and surprising happens in higher dimensions---has the flavor of real mathematics. But once again, the methods used to settle the last case put the final result squarely in the realm of computer science. Propositional satisfiability solvers are used to solve problems in a wide range of domains, including hardware and software verification as well as discrete optimization. Mathematicians do not conventionally view search problems as being a part of mathematics, especially not when the state space is finite. One mathematician I know, upon hearing of the resolution of Keller's conjecture, said to me: ``don't expect mathematicians to be interested---this is a `finite computation.'\,''

Theorem~\ref{theorem:hales}, the Kepler conjecture, rises higher in the ranking. After all, it isn't easy to ignore the proof a conjecture that remained open for nearly four centuries. The proof did not involve heuristic search, but it did involve lengthy calculations, including brute force enumeration and numerical bounds on large numbers of nonlinear expressions. That is enough to convince many mathematics that the problem is not really a mathematical one, but rather, a matter of engineering. Appel and Haken's proof of the four color theorem was met with a similar response, highlighting the concern that computational methods do not deliver a suitable mathematical understanding.

In 2018, I wrote an article, \emph{The Mechanization of Mathematics} \cite{avigad:18}, that surveyed recent logic-based uses of computers in mathematics. A colleague commented on a draft:
\begin{quote}
\ldots the most striking point is how your emphasis is absolutely orthogonal to what interests me (and other mathematicians \ldots).

Let me explain this, exaggerating slightly for comparison. First, for the purpose of this discussion, what I call a traditional mathematician is someone who has a permanent position that involves proving theorems
\ldots\,. What they do is try to get a conceptual and rigorous understanding of which mathematical statements are true. Both adjectives are important. To them, Hales' proof of the Kepler conjecture is nothing like solving the conjecture. And the Flyspeck project is purely computer science. That doesn't mean it is not interesting; it is something different, because at least part of the proof lacks the ``conceptual'' adjective. You will often hear mathematicians, talking about their own proofs, saying things like ``there is nothing to understand here, it's only a computation.''

To most mathematicians, the interesting part of a proof is the high level structure, and the technical details are only a chore that is necessary to make sure we haven't missed anything important.
\end{quote}
Notice that the passage emphasizes that the attitude it describes is not a comment on how interesting and important the proof of the Kepler conjecture is, but, rather, the extent to which it belongs to traditional mathematics.

Viazovska's result on the optimality of $E_8$ and the follow-up result on the Leech lattice fare the best under the standards of Section~\ref{section:understanding}. There are a number of strikes against them: the proofs are short, they involve a fair amount of calculation, they do not introduce complicated new definitions, and the discoveries build on intuitions gleaned from numerical computation and experimentation. So what makes them real mathematics? The fact that $E_8$ and the Leech lattice are richly complicated and highly symmetric objects, and the fact that they many have interesting algebraic properties, suggest that they, themselves, are conceptually deep. It also helps that the proof draws on modular forms, representation theory, and harmonic analysis, themselves well established as central components of traditional mathematics. Cohn's description of Viazovska's work emphasizes these features:

\begin{quote}

\smallskip

``\ldots it's wonderful to see a relatively simple proof of a deep theorem in sphere packing.''

\smallskip

``Her proof is thus a notable contribution to the story of $E_8$, and more generally the story of exceptional structures in mathematics.''

\smallskip

``Viazovska \ldots establishes a new connection between modular forms and discrete geometry.''

\smallskip

``Instead of justifying sphere packing by aspects of the problem or its applications, we'll justify it by its solutions: a question is good if it has good answers. Sphere packing turns out to be a far richer and more beautiful topic than the bare problem statement suggests. From this perspective, the point of the subject is the remarkable structures that arise as dense sphere packings.''

\smallskip

``\ldots modular forms are deep and mysterious functions connected with lattices, as are the magic functions, so wouldn't it make sense for them to be related?''

\smallskip

``Despite our lack of understanding [of other higher dimensional sphere packing problems], the special role of eight and twenty-four dimensions aligns with our experience elsewhere in mathematics. Mathematics is full of exceptional or sporadic phenomena that occur in only finitely many cases, and the $E_8$ and Leech lattices are prototypical examples. These objects do not occur in isolation, but rather in constellations of remarkable structures.''

\smallskip

\end{quote}
In sum, what makes the results attractive are the structures involved and the rich connections to other parts of mathematics. Even the title of the article, ``A conceptual breakthrough in sphere packing,'' emphasizes the conceptual aspect.

\section{Varieties of understanding}
\label{section:discussion}

In Section~\ref{section:understanding}, I characterized a view of mathematical understanding, and in Section~\ref{section:assessment}, I observed that the results of Section~\ref{section:theorems} do not rank particularly well by that standard. What are we to make of this?

Let's clarify what is at stake. At issue is not the quality or importance of the results; the theorems are impressive, and mathematicians are generally happy to grant them the respect they deserve. The questions are natural, the theorems are pretty, the methods are clever, and the insights are substantial. The mathematical community is open-minded and generous, and it is quick to appreciate the value of mathematical results in physics, chemistry, computer science, and even discrete geometry, whether or not they are seen as core mathematics.

Nor is it a question of popularity. The results of Section~\ref{section:theorems} have all made it into the popular press and have gained a level of recognition that very few core mathematical results ever receive. It is not a matter of institutional support; most of the authors cited in Section~\ref{section:theorems} have good jobs and stable academic positions. The question is not about academic recognition; the results of Section~\ref{section:theorems} have high citation counts, generally much higher than papers associated with, say, the Langlands problem. Nor is it a question of funding; grants are generally larger and easier to come by in applied mathematics, computer science, and engineering, and availability of funding in areas of pure mathematics tends to track their proximity to the more applied branches. By that standard, research in discrete geometry does just fine.

Moreover, I don't see the mathematical elite, or anyone making pronouncements as to what traditional mathematicians do, complaining about any of that. Mathematicians are often proud to distance themselves from the claim that anything they do should be considered practically useful. People don't go into pure mathematics because they are looking for fame and fortune. They do it because they find it enjoyable, meaningful, and satisfying; because they admire and respect the power of mathematical thought; because they see themselves as part of a long and worthy tradition; and because they find the subject elegant and beautiful.

The question I am asking here is therefore more of a philosophical one. Insofar as there is an essence of mathematical understanding, is it adequately captured by the conceptual view described in Section~\ref{section:understanding}? Or would a better articulation of this essence leave more room for the results of Section~\ref{section:theorems}? Nobody denies that the results are rigorous, and nobody denies that they introduce interesting new ideas. They certainly meet some of the criteria for good mathematics: they open up new capacities for rigorous thought; they provide new means of solving hard problems; they offer new representations, conceptualizations, and methods; and they raise new questions and avenues for exploration.

Let me replace the broad philosophical question with two that are more focused:
\begin{itemize}
  \item Can the kinds of computational methods used to obtain the results of Section~\ref{section:theorems} contribute to core areas of mathematics and help answer some of the questions that are currently viewed as important?
  \item Can/should/will computational methods like the ones used toward obtaining the results of Section~\ref{section:theorems} become, in and of themselves, part of a proper mathematical understanding?
\end{itemize}

The answer to the first question is likely to be ``yes.'' Computers fundamentally change the kinds of reasoning we can perform and dramatically increase the level of complexity of the inferences that we can carry out with any measure of reliability. The Langlands program, for example, seeks correspondences between certain types of number-theoretic and geometric data. Perhaps there is an aesthetic judgment as to whether a proposed correspondence is sufficiently natural, but insofar as the problem can be made precise, I see no reason to think that search techniques can't be used to supplement conventional mathematical insight.

This last claim may seem ludicrous to practitioners. Computers are certainly nowhere near being able to construct the kinds of conceptual definitions and proofs that experts do. But let's recognize that we have been looking for those particular kinds of definitions and proofs because we \emph{can}; they are the particular means we have developed to meet mathematical challenges with the cognitive resources at our disposal. The computer is a new type of resource, and it can produce different kinds of solutions. Depending on how you count, we have been using the conceptual approach for a century or more. In contrast, we still have a lot to learn about how computers can be used to address core mathematical problems. Maybe it will require a modern-day Gauss to figure out how to harness the power of computation, or maybe we will get there in incremental steps. But it is impossible to deny that raw computing power extends our mathematical capacities in a significant way, and it is only a matter of time before we learn how to make proper use of it.

Whereas the first question is pragmatic, the second is a matter of mathematical values. If it turns out that some Langlands-like questions can be answered with the use of computation, there is always the possibility that the mathematical community will interpret this as a demonstration that, in hindsight, the Langlands program is not as deep as we thought it was. There is always room to say, ``Aha! Now we see that it is just a matter of computation.'' The question is whether computational methods might one day become \emph{fashionable} in Buzzard's sense. If computational methods prove useful in solving problems that were antecedently deemed important, questions about the new methods themselves might be judged to be important in their own right. Can we imagine that happening?

The history of mathematics offers some useful lessons in that regard. The first lesson is that mathematical values change over time, not just views as to what counts as rigorous, but also views as to what counts as a proper understanding. Consider, for example, the fact that in the Western tradition, from ancient times to the seventeenth century, geometry was held to be the proper foundation for mathematics. Today, we think of an equation like $y = x^2 + 2 x + 1$ as expressing a relationship between elements of the real number structure, a relationship that can be depicted geometrically by graphing it. For the early algebraists, the perspective was inverted: variables stood for geometric magnitudes like lengths and areas, and algebraic manipulations were abstract proxies for geometric construction. In 1637, Descartes' landmark \emph{G\'eom\'etrie} opened up new opportunities for the use of algebraic methods in solving geometric problems, but even for Descartes, such a problem was not solved until the algebraic solution was reinterpreted as a geometric construction.

Such a detour through algebraic methods was unacceptable to Newton, who felt that although such methods could be used as a heuristic for finding solutions, a proposition was not suitably demonstrated, and not properly understood, until it was given a purely geometric presentation. It is well known that Newton developed an algebraic version of calculus around 1665, but suppressed publication in favor of the more geometric methods of his monumental \emph{Philosophi{\ae} Naturalis Principia Mathematica}, published in 1687. Newton's antipathy to algebraic methods has been well documented by Guicciardini \cite{guicciardini:99,guicciardini:09}. For example, in Newton's writings of the 1690s, we find:
\begin{quote}
Men of recent times, eager to add to the discoveries of the ancients, have united specious arithmetic with geometry. Benefitting from that, progress has been broad and far reaching if your eye is on the profuseness of output but the advance is less of a blessing if you look at the complexity of its conclusions. For these computations, progressing by means of arithmetical operations alone, very often express in an intolerably roundabout way quantities which are designated by the drawing of single line.\footnote{This passage, from Newton's \emph{Mathematical Papers} \cite[vol.~7, p.~251]{newton:mp}, appears in Guicciardini \cite[Section 4.6]{guicciardini:09}, along with others of a similar character. See also \cite[Chapter 2]{guicciardini:09}.}
\end{quote}
He stated such a position even more strongly in the mid-1670s, when he judged the ancient approach to solving a problem by Pappus to be ``more elegant by far than the Cartesian one,'' since Descartes
\begin{quote}
achieved the results by an algebraic calculus which, when transposed into words (following the practice of the Ancients in their writings) would prove to be so tedious and entangled as to produce nausea, nor might it be understood.\footnote{From \cite[vol.~4, p.~277]{newton:mp}, translated by Westfall \cite[p.~379]{westfall:80}, and quoted by Guicciardini \cite[p.~29]{guicciardini:99}}
\end{quote}
In late 1683 or early 1684, Newton lectured on the ``treatment of composition,'' which Guicciardini describes as ``the synthetic, constructive phase of the problem-solving process.'' According to Newton:
\begin{quote}
Analysis guides us to the composition, but true composition is not achieved before it is freed from analysis. Let even the slightest trace of analysis be present in the composition and you will not yet have attained true composition.\footnote{From \cite[vol.~4, p.~477]{newton:mp}, quoted by Guicciardini \cite[pp.74--75]{guicciardini:09}.}
\end{quote}

As Guicciardini's scholarship shows, Newton's views were based in part on somewhat mystical attributions to the wisdom of the ancients, but they were also motivated by deep-seated views about mathematical meaning. Today, we tend to view the translation of a visual argument into set-theoretic terms as a step toward rigor. From Newton's standpoint, this is a step \emph{away} from a proper mathematical understanding. It may be tempting to write Newton off as more of a physicist than a mathematician, but that is not easy to do: throughout his career he held mathematics to be prior to its physical applications, and his insistence on rigor and purity of method, as well as his ability to solve fiendishly difficult mathematical problems, mark him as one of our own.

The nineteenth-century transition to the modern conceptual view was beset by similar misgivings. The conceptual approach was designed to address situations where calculational methods were straining to reign in complexity in the theory of equations, in analysis, and in number theory. In analogy to the seventeenth-century view of the use of algebraic methods in geometry, the new abstract methods were seen as heuristic shortcuts to obtaining results that were properly justified by explicit calculation. Once again, the community was initially reluctant to stray too far. Today it is hard to appreciate the extent to which calculation was central to the lives of nineteenth-century mathematicians, even the conceptual innovators. Gauss famously recounted, in a letter to his student, Johann Franz Encke, the great pleasure it gave him to spend an idle quarter hour tabulating frequencies of prime numbers.\footnote{See, for example, \cite{tschinkel:06}.} Riemann's short paper on the zeta function was a tour-de-force of calculation, and Dedekind motivated his algebraic approaches to number theory by working out concrete examples. Calculation has always been central to mathematics; we can perhaps forgive nineteenth-century mathematicians for thinking that symbolic representations and algorithms are an essential part of mathematical understanding.

A second lesson from history is that when mathematics does change, preconceived notions as to what constitutes understanding are invariably pushed aside by mathematical success. Mathematics isn't beholden to dogma; we may think we know what constitutes good mathematics, but mathematics may tell us otherwise. In the eighteenth century, the algebraic approach to calculus gave us Euler, the Bernoullis, Lagrange, and Laplace, and it gave us advances in mechanics, probability, and physics. It is no wonder we were willing to turn our backs on geometry. The abstract, infinitary methods of the nineteenth century gave us algebraic number fields, complex analysis, manifolds, Lie groups, and the beginnings of algebraic geometry. Concerns about computational meaning simply took a back seat.

It is perhaps the fate of mathematicians through the ages to think that their views of mathematics are timeless. But today the notion that mathematical truth is grounded in geometry seems quaint, and contemporary mathematicians scoff at the suggestion that the subject is fundamentally about calculation. That should remind us that mathematicians of the future may one day judge our present view that the essence of mathematics lies in a certain type of conceptual understanding---built on networks of complicated algebraic definitions, locally surveyable, with communities of experts straining to hold the pieces together---as equally quaint and outdated. To be sure, the conceptual view has won us great success in answering traditional questions and generating new ones. But as new methods arise and fashions change, our views of mathematics will change as well.

And so our current preoccupation with conceptual methods could well give way to a more expansive view of mathematical understanding, one in which the computer plays a more central role. This isn't to say that we will then abandon the conceptual viewpoint entirely. A third lesson from the history of mathematics is that, as the subject evolves, it has a prodigious capacity for reinterpreting its past. We still recognize the importance of Euclid's \emph{Elements}, Archimedes' calculations of volumes, and Appolonius' work on conics; we acknowledge the brilliance of Newton and Leibniz in the uniform methods they developed to calculate areas and tangents; we appreciate the delicate arguments used by Abel and Galois to demonstrate the unsolvability of the quintic; and we value Gauss' detailed manipulations of quadratic forms. We see all of these as part of the modern, conceptual understanding, or, at least, as an essential part of its history. If eventually we end up with a view of mathematical understanding that can better accommodate the use of computation, it will likely be seen as an inevitable refinement of the conceptual point of view.

\section{Conclusions}

By highlighting the fact that views of mathematical understanding change over time, I do not mean to deny that mathematics has any fixed characteristics. The longevity and continuity of the subject, as well as our ability to understand, appreciate, and build on the results of the past, show that there is something about mathematics that manages to transcend the vicissitudes of time, place, and culture. It has always been concerned with solving hard problems, communicating ideas in precise ways, and establishing stable foundations for scientific reasoning. The subject aims to establish coordination between practitioners not only at a given point in time, but across centuries. Its remarkable success at doing so is evidence that mathematics draws on something fundamental about the way we cognize the world around us.

What changes, rather, are our judgments as to the specific ways we can best achieve our mathematical goals. These judgments depend on historical circumstances, including the body of mathematical knowledge accumulated to date, the demands of science, social and institutional factors, and, yes, technology. The only sure-fire way of finding out where mathematics is going is to wait and see, but we ourselves are part of the mathematical tradition, and our philosophical views play a role in guiding its development. Being mindful of these views and discussing them is an important part of trying to understand what constitutes good mathematics.

Potential uses of computational methods, in particular, should challenge us to think about our mathematical goals and how they are best served. It is notable that the headline photograph to Cohn's \emph{Notices} article features Viazovska at her computer, immersed in thought, with one hand on a mouse and with a coffee cup nearby. The caption reads ``Maryna Viazovska solved the sphere packing problem in eight dimensions.'' For all we know, in the picture, she may be catching up on email or surfing the web for sports results. The caption, however, suggests that she is using the computer to unlock deep mathematical secrets. The iconic trope of a mathematician standing before a blackboard covered with symbols and arrows may be giving way to the image of a mathematician at a keyboard, coaxing mathematical understanding from calculation, simulation, and search. That should give us something to think about.

\bibliographystyle{amsplain}
\bibliography{varieties_arxiv}

\providecommand{\bysame}{\leavevmode\hbox to3em{\hrulefill}\thinspace}
\providecommand{\MR}{\relax\ifhmode\unskip\space\fi MR }
\providecommand{\MRhref}[2]{%
  \href{http://www.ams.org/mathscinet-getitem?mr=#1}{#2}
}
\providecommand{\href}[2]{#2}
\begin{thebibliography}{10}

\bibitem{avigad:06}
Jeremy Avigad, \emph{Methodology and metaphysics in the development of
  {D}edekind's theory of ideals}, The Architecture of Modern Mathematics
  (Jos\'e Ferreir\'os and Jeremy Gray, eds.), Oxford Univ. Press, Oxford, 2006,
  pp.~159--186. \MR{2258019}

\bibitem{avigad:18}
\bysame, \emph{The mechanization of mathematics}, Notices Amer. Math. Soc.
  \textbf{65} (2018), no.~6, 681--690. \MR{3792862}

\bibitem{biere:et:al:09}
Armin Biere, Marijn Heule, Hans van Maaren, and Toby Walsh (eds.),
  \emph{Handbook of satisfiability}, {IOS} Press, 2009.

\bibitem{borwein:bailey:08}
Jonathan Borwein and David Bailey, \emph{Mathematics by experiment}, second
  ed., A K Peters, Ltd., Wellesley, MA, 2008, Plausible reasoning in the 21st
  Century. \MR{2473161}

\bibitem{borwein:devlin:09}
Jonathan Borwein and Keith Devlin, \emph{The computer as crucible: An
  introduction to experimental mathematics}, A K Peters, Ltd., Wellesley, MA,
  2009. \MR{2464847}

\bibitem{brakensiek:et:al:20}
Joshua Brakensiek, Marijn Heule, John Mackey, and David Narv{\'{a}}ez,
  \emph{The resolution of {Keller's} conjecture}, International Joint
  Conference on Automated Reasoning (IJCAR) 2020 (Nicolas Peltier and Viorica
  Sofronie{-}Stokkermans, eds.), Springer, 2020, pp.~48--65.

\bibitem{chen:engel:glotzer:10}
Elizabeth~R. Chen, Michael Engel, and Sharon~C. Glotzer, \emph{Dense
  crystalline dimer packings of regular tetrahedra}, Discrete Comput. Geom.
  \textbf{44} (2010), no.~2, 253--280. \MR{2671012}

\bibitem{cohn:17}
Henry Cohn, \emph{A conceptual breakthrough in sphere packing}, Notices Amer.
  Math. Soc. \textbf{64} (2017), no.~2, 102--115. \MR{3587715}

\bibitem{cohn:et:al:17}
Henry Cohn, Abhinav Kumar, Stephen~D. Miller, Danylo Radchenko, and Maryna
  Viazovska, \emph{The sphere packing problem in dimension 24}, Ann. of Math.
  (2) \textbf{185} (2017), no.~3, 1017--1033. \MR{3664817}

\bibitem{conway:torquato:06}
J.~H. Conway and S.~Torquato, \emph{Packing, tiling, and covering with
  tetrahedra}, Proc. Natl. Acad. Sci. USA \textbf{103} (2006), no.~28,
  10612--10617. \MR{2242647}

\bibitem{gauss:disquisitiones:arithmeticae}
Carl~Friedrich Gauss, \emph{Disquisitiones arithmeticae}, Yale University
  Press, New Haven, Conn.-London, 1966, Translated into English by Arthur A.
  Clarke. \MR{0197380}

\bibitem{gowers:00}
W.~T. Gowers, \emph{The two cultures of mathematics}, Mathematics: frontiers
  and perspectives, Amer. Math. Soc., Providence, RI, 2000, pp.~65--78.
  \MR{1754768}

\bibitem{gowers:10}
\bysame, \emph{Decompositions, approximate structure, transference, and the
  {H}ahn-{B}anach theorem}, Bull. Lond. Math. Soc. \textbf{42} (2010), no.~4,
  573--606. \MR{2669681}

\bibitem{green:tao:08}
Ben Green and Terence Tao, \emph{The primes contain arbitrarily long arithmetic
  progressions}, Ann. of Math. (2) \textbf{167} (2008), no.~2, 481--547.
  \MR{2415379}

\bibitem{guicciardini:99}
Niccol\`o Guicciardini, \emph{Reading the \emph{{P}rincipia}: The debate on
  {N}ewton's mathematical methods for natural philosophy from 1687 to 1736},
  Cambridge University Press, Cambridge, 1999. \MR{1725356}

\bibitem{guicciardini:09}
\bysame, \emph{Isaac {N}ewton on mathematical certainty and method}, MIT Press,
  Cambridge, MA, 2009. \MR{2562084}

\bibitem{haji:akbari:et:al:09}
Amir Haji-Akbari, Michael Engel, Aaron~S. Keys, Xiaoyu Zheng, Rolfe~G.
  Petschek, Peter Palffy-Muhoray, and Sharon~C. Glotzer, \emph{Disordered,
  quasicrystalline and crystalline phases of densely packed tetrahedra}, Nature
  \textbf{462} (2009), 773--777.

\bibitem{hales:et:al:17}
Thomas Hales, Mark Adams, Gertrud Bauer, Tat~Dat Dang, John Harrison, Le~Truong
  Hoang, Cezary Kaliszyk, Victor Magron, Sean McLaughlin, Tat~Thang Nguyen,
  Quang~Truong Nguyen, Tobias Nipkow, Steven Obua, Joseph Pleso, Jason Rute,
  Alexey Solovyev, Thi Hoai~An Ta, Nam~Trung Tran, Thi~Diep Trieu, Josef Urban,
  Ky~Vu, and Roland Zumkeller, \emph{A formal proof of the {K}epler
  conjecture}, Forum Math. Pi \textbf{5} (2017), e2, 29. \MR{3659768}

\bibitem{hales:05}
Thomas~C. Hales, \emph{A proof of the {K}epler conjecture}, Ann. of Math. (2)
  \textbf{162} (2005), no.~3, 1065--1185. \MR{2179728}

\bibitem{hamami:morris:20}
Yacin Hamami and Rebecca~L. Morris, \emph{Philosophy of mathematical practice:
  a primer for mathematics educators}, ZDM Mathematics Education \textbf{52}
  (2020), 1113--1126.

\bibitem{keller:30}
Ott-Heinrich Keller, \emph{\"{U}ber die l\"{u}ckenlose {E}rf\"{u}llung des
  {R}aumes mit {W}\"{u}rfeln}, J. Reine Angew. Math. \textbf{163} (1930),
  231--248. \MR{1581241}

\bibitem{kisielewicz:17}
Andrzej~P. Kisielewicz, \emph{Rigid polyboxes and {K}eller's conjecture}, Adv.
  Geom. \textbf{17} (2017), no.~2, 203--230. \MR{3652241}

\bibitem{lagarias:shor:92}
Jeffrey~C. Lagarias and Peter~W. Shor, \emph{Keller's cube-tiling conjecture is
  false in high dimensions}, Bull. Amer. Math. Soc. (N.S.) \textbf{27} (1992),
  no.~2, 279--283. \MR{1155280}

\bibitem{mackey:02}
John Mackey, \emph{A cube tiling of dimension eight with no facesharing},
  Discrete Comput. Geom. \textbf{28} (2002), no.~2, 275--279. \MR{1920144}

\bibitem{perron:40a}
Oskar Perron, \emph{\"{U}ber l\"{u}ckenlose {A}usf\"{u}llung des
  {$n$}-dimensionalen {R}aumes durch kongruente {W}\"{u}rfel}, Math. Z.
  \textbf{46} (1940), 1--26. \MR{3041}

\bibitem{perron:40b}
\bysame, \emph{\"{U}ber l\"{u}ckenlose {A}usf\"{u}llung des {$n$}-dimensionalen
  {R}aumes durch kongruente {W}\"{u}rfel. {II}}, Math. Z. \textbf{46} (1940),
  161--180. \MR{2185}

\bibitem{szabo:86}
S.~Szab\'{o}, \emph{A reduction of {K}eller's conjecture}, Period. Math.
  Hungar. \textbf{17} (1986), no.~4, 265--277. \MR{866636}

\bibitem{tschinkel:06}
Yuri Tschinkel, \emph{About the cover: on the distribution of primes---{G}auss'
  tables}, Bull. Amer. Math. Soc. (N.S.) \textbf{43} (2006), no.~1, 89--91.
  \MR{2201552}

\bibitem{viazovska:17}
Maryna~S. Viazovska, \emph{The sphere packing problem in dimension 8}, Ann. of
  Math. (2) \textbf{185} (2017), no.~3, 991--1015. \MR{3664816}

\bibitem{westfall:80}
Richard~S. Westfall, \emph{Never at rest: A biography of {I}saac {N}ewton},
  Cambridge University Press, Cambridge, 1980. \MR{741027}

\bibitem{newton:mp}
D.~T. Whiteside (ed.), \emph{The mathematical papers of {I}saac {N}ewton},
  Cambridge University Press, Cambridge, 1967--1981, in seven volumes.
  \MR{0505130}

\end{thebibliography}

\end{document}